\documentclass[12pt]{amsart} 
\usepackage{amssymb,amsmath,amscd}  
 
\newtheorem{thm}{Theorem}[section] 
\newtheorem{corollary}[thm]{Corollary} 
\newtheorem{prop}[thm]{Proposition} 
\newtheorem{lemma}[thm]{Lemma} 
\newtheorem{fact}[thm]{Fact} 
 
\theoremstyle{definition} 
\newtheorem{defn}[thm]{Definition} 
\newtheorem{example}[thm]{Example} 
 
\theoremstyle{remark} 
\newtheorem{remark}[thm]{Remark} 
 
\newcommand{\bt}{\begin{thm}} 
\newcommand{\et}{\end{thm}} 
\newcommand{\bp}{\begin{prop}} 
\newcommand{\ep}{\end{prop}} 
\newcommand{\bd}{\begin{defn}} 
\newcommand{\ed}{\end{defn}} 
\newcommand{\bl}{\begin{lemma}} 
\newcommand{\el}{\end{lemma}} 
\newcommand{\bfa}{\begin{fact}} 
\newcommand{\efa}{\end{fact}} 
\newcommand{\bc}{\begin{corollary}} 
\newcommand{\ec}{\end{corollary}} 
\newcommand{\bex}{\begin{example}} 
\newcommand{\eex}{\end{example}} 
\newcommand{\br}{\begin{remark}} 
\newcommand{\er}{\end{remark}}

\newcommand{\isom}{\cong} 
\newcommand{\dual}[1]{{#1}^\vee} 
 
\newcommand{\Hom}{{\mathcal H}om}

\newcommand{\sotto}[2]{#1_{#2}}  
 
\newcommand{\rrr}{\rightarrow}

\newcommand{\ideal}[1]{\sotto {{\mathcal I}}{#1}}

\newcommand{\exact}[3] 
{0 \rrr #1 \rrr #2 
\rrr #3 \rrr 0}

\newcommand{\pso}{\mathbb{P}^3}

\newcommand{\ptwo}{\mathbb{P}^2}

\newcommand{\Z}{\mathbb{Z}}

\newcommand{\ccc}{{\mathcal C}} 
 
\newcommand{\coo}{{\mathcal O}}

\newcommand{\caf}{{\mathcal F}}

\newcommand{\scrl}{{\mathcal L}} 
\newcommand{\cae}{{\mathcal E}}

\newcommand{\df}{\partial} 
 
\newcommand{\gmc}{\gamma_{C}} 
\newcommand{\gmz}{\gamma_{Z}}

\begin{document}

\title{Curves in the double plane} 
 
\author{Robin Hartshorne} 
\address{Department of Mathematics, University of California,  
Berkeley CA94720, USA} 
\email{robin@math.berkeley.edu} 
 
\author{Enrico Schlesinger} 
\address{Dipartimento di Matematica, Universit\`{a} degli Studi di Trento, Via Sommarive 14, 38050 POVO (Trento), Italia} 
\curraddr{Dipartimento di Matematica, Politecnico di Milano, Via Bonardi 9, 20133 Milano, Italia} 
\email{schlesin@science.unitn.it,enrsch@mate.polimi.it} 
\thanks{The second author was supported by a post-doctoral  
fellowship of the Dipartimento di Matematica, Universit\`{a} degli Studi 
di Trento} 
 

\begin{abstract} 
We study in detail locally  
Cohen-Macaulay curves in ${\mathbb{P}^3}$ 
which are contained in a double plane $2H$, thus 
completing the classification of curves lying on surfaces  
of degree two. We describe the irreducible components 
of the Hilbert schemes $H_{d,g} (2H)$ of locally  
Cohen-Macaulay curves in $2H$ of degree $d$ and arithmetic genus $g$, 
and we show that $H_{d,g} (2H)$ is connected.    
We also discuss the Rao module of these curves and  
liaison and biliaison equivalence classes.  
\end{abstract} 
 
\maketitle 
\section{Introduction} \label{one}  
Much attention has been given in recent years to the classification 
of curves in projective space. Here we define a curve to be a  
purely one-dimensional locally Cohen-Macaulay (i.e. without  
embedded points) closed subscheme of $ \pso_k$, the projective  
three-space over an algebraically closed field $k$. Our goal in  
this paper is to answer all the interesting questions about curves  
contained in surfaces of degree two. We succeed, with one notable 
exception.  
 
Curves in the nonsingular quadric  surface $Q$ and in the quadric cone 
$Q_0$ are well-known  \cite{AG} (Exercises III 5.6 and  V 2.9). Curves 
in the union of two planes  $H_1 \cup H_2$ were studied in \cite{hgd}, 
section  5.  The  main contribution  of this  paper is  the systematic 
study of curves contained in a double plane $2H$. 
 
We are also interested in flat families of curves, so that we can  
identify the irreducible components of the Hilbert scheme. Families 
of curves whose general member lies on a nonsingular quadric $Q$ 
and whose special member lies on a cone $Q_0$ or the union of  
two planes $H_1 \cup H_2$ were studied in \cite{zeuthen}. In this paper  
we will identify the irreducible components of the scheme  
$H_{d,g} (2H)$ of curves of degree $d$ and genus $g$ lying in a  
double plane, and specializations between them. The question we have  
not been able to answer is what families are there whose general  
member lies on a nonsingular quadric and whose special member lies in 
$2H$? 
For example, we do not know if there is a family specializing from  
four skew lines on a quadric to a curve contained in the double plane.  
 
Now let us describe the contents of the current paper. To each  
curve $C$ contained in the double plane $2H$, we assign a triple  
$\{Z,Y,P\}$ consisting of a zero-dimensional closed subscheme $Z$  
in the (reduced) plane $H$, and two curves $Y,P$ in $H$. Roughly  
speaking, $P$ with embedded points $Z$ is the intersection of $C$ with  
$H$, 
while $Y$ is the residual intersection (Proposition 2.1).  
Conversely to each such triple, satisfying certain conditions, we  
can associate curves $C$ in $2H$ (Proposition 3.1). The numerical  
invariants of $C$ can be computed in terms of $Z,Y,P $ (section~\ref{five}). 
Thus the triple $\{Z,Y,P \}$ of data in the plane $H$ is  
an effective tool for studying curves $C$ in the double plane.  
 
We can translate information about a flat family of  
curves in $2H$ to families of triples in $H$ (section \ref{four}). This  
allows us to identify the irreducible components of the Hilbert  
scheme $H_{d,g} (2H)$. They are given by triples of integers  
$(z,y,p)$ representing the length of $Z$ and the degrees of  
$Y$ and $P$ (Theorem~\ref{components}).  
 
We also discuss the Rao module of these curves (section \ref{five}) and  
liaison and biliaison equivalence classes (section \ref{six}). Finally,  
generalizing a technique of Nollet \cite{nthree}, we show the existence of  
flat families joining the various irreducible components of  
$H_{d,g} (2H)$, and conclude that $H_{d,g} (2H)$ is connected 
(Theorem~\ref{connect}).  
 
Our purpose in studying curves in the double plane  
was to treat one very special case in view of the more  
general problem whether the Hilbert scheme $H_{d,g}$ of all 
locally Cohen-Macaulay curves in $\pso$ is connected for  
all $d$,$g$. Because of our connectedness result (\ref{connect}), it  
would be sufficient for that more general problem to show that any  
curve in $\pso$ can be connected by a sequence of flat families to a 
curve in 2H.  
 
We would like to thank E. Ballico for encouraging us to write this paper.

\section{The triple associated to a curve in 2H} \label{two} 
We first recall the notion  of ``residual scheme''  
(\cite{fulton}, 9.2.8, and \cite{PS} ). Suppose $T \subseteq W$ 
are closed subschemes of an ambient scheme $V$. The residual scheme 
$R$ of $T$ in $W$ is the closed subscheme of $V$ with ideal sheaf  
$$ 
\ideal{R} = (\ideal{W} : \ideal{T}). 
$$ 
Intuitively, $R$ is equal to $W$ minus $T$. It does not depend on the 
ambient scheme. We will need this notion in the following two cases. 
If $T$ is a Cartier divisor in $V$, then 
$\ideal{R} \ideal{T} = \ideal{W}$. On the other hand, suppose 
that $H$ is an effective Cartier divisor on $V$, and $T$ is the scheme 
theoretic intersection of $W$ and $H$. In this case, we say  
(cf. \cite{hhdg}, p.176) that  $R$ is the residual scheme to the 
intersection of  $W$  with  $H$. 
We have 
$$ 
\ideal{R} \ideal{H} = \ideal{W} \cap \ideal{H} 
$$ 
hence an exact sequence: 
\begin{equation} \label{res-int} 
 0 \rrr  
\ideal{R}(-H) \rrr \ideal{W} \rrr  
\ideal{W \cap H,H}  \rrr 0.  
\end{equation} 
 
Let $H$ be a plane in $\pso_{k}$, 
defined by $h=0$.  Let $C$ be a  
curve contained in  the scheme $2H$ (``the double  plane'') defined by 
the equation  $h^2 = 0$. To the  curve $C$ we will  associate a triple 
$T(C) = \{Z,Y,P \}$ where  $Z$ is a zero-dimensional subscheme of $H$, 
and $Y$,$P$ are curves in $H$, with $Z \subseteq Y \subseteq P$. 
 
First consider the scheme-theoretic intersection $C \cap H$. This will 
be a one dimensional subscheme of  $H$, possibly with embedded points.  
So we can write 
$$ 
\ideal{C \cap H,H} = \ideal{Z,H} (-P) 
$$ 
where $Z$ is zero-dimensional and $P$ is a curve.  
In fact, the inclusion 
$\ideal{C \cap H,H} \hookrightarrow  \coo_{H}$ defines a  global 
section of the invertible sheaf $\Hom (\ideal{C \cap H,H}, \coo_{H})$ 
whose scheme of zeros is the effective Cartier divisor $P \subset H$, 
and $Z$ is the residual scheme of $P$ in $C \cap H$. 
 
Next  we  let  $Y$  be  the  residual scheme to the 
intersection  of  $C$  with $H$:  
it is a  curve in $\pso$. By the discussion above we have an exact sequence 
\begin{equation} \label{first} 
 0 \rrr  
\ideal{Y,\pso}(-1) \stackrel{h}{\rrr} \ideal{C,\pso} \rrr  
\ideal{C \cap H,H}  \rrr 0.  
\end{equation} 
Since  $C$   is  contained   in  $2H$,  $Y$   will  be   contained  in 
$H$.  Furthermore, since $P$  is  the  largest curve  in  $H$ which  is 
contained in $C$, it is clear that $Y \subseteq P$. Note that $C$ is  
contained in the reduced plane $H$ if and only if $Y$ is empty. 
 
Now we use the inclusion $P \subseteq C \cap H$  to create a diagram 
as follows   
\begin{equation} \label{second} 
\begin{CD} 
&&0 && 0 && 0&& \\ 
&&@VVV   @VVV @VVV \\ 
0 @>>> \ideal{Y,\pso} (-1)     @>h>>    \ideal{C, \pso}  @>>>   
\ideal{C \cap H,H}   @>>>  0  \\ 
&&@VVV   @VVV @VVV \\ 
0 @>>> \coo_{\pso} (-1)  @>h>>  \ideal{P,\pso} @>>>  \ideal{P,H} @>>> 0\\ 
&&@VVV   @VV{u}V  @VVV \\ 
0 @>>> \coo_{Y} (-1) @>>> \scrl @>>>  \coo_{Z} (-p) @>>> 0    \\ 
&&@VVV  @VVV @VVV \\ 
&&0 && 0 && 0&& \\ 
\end{CD} 
\end{equation} 
where $p$ is the degree of $P$. 
Since $C$ is locally Cohen-Macaulay, the ideal sheaf $\ideal{C,\pso}$ 
has depth $\geq 2$ at every closed point. Therefore the sheaf $\scrl$, 
defined as the quotient $\ideal{P} / \ideal{C}$, will have depth $\geq 
1$. Multiplying a section of $\scrl$ by a local section of $\ideal{Y}$ 
will give something zero outside of $Z$. But $\scrl$ can have no section 
supported at a point because of its depth. So we see that $\scrl$ is an 
$\coo_{Y}$-module. Because $Y$ is a Gorenstein curve, $\scrl$ is a 
reflexive sheaf of rank one \cite{hgd}, 1.6, and then by \cite{hgd}, 2.8 and 2.9,  
$\scrl$ is 
of the form $\scrl(D-1)$ for some effective generalized divisor $D$ on 
$Y$. Then by  \cite{hgd}, 2.10, the quotient of $\scrl$ by $\coo_{Y} 
(-1)$ is 
$$ 
\coo_{Z} (-p) \cong \omega_{D} \otimes \dual{\omega}_{Y} (-1), 
$$ 
so 
$$ 
\coo_{Z} (-p) \cong \omega_{D} (2-y) 
$$ 
where $y$ is the degree of $Y$. 
 
From  this  we  conclude  several  things. First  $Z=D$,  so  that  $Z 
\subseteq   Y$.   Secondly,  $\coo_{Z}$   is  locally   isomorphic  to 
$\omega_{Z}$,  so  $Z$  is  a  Gorenstein scheme.   Since  $Z$  is  of 
codimension two  in $H$,  it follows from the theorem of Serre 
\cite{eisenbud}, Corollary  21.20,  that 
$Z$ is a locally complete intersection in $H$. Summing up,  
\bp \label{1.1} 
To each curve $C$  in $2H$, we 
can  associate a  triple  $T(C) =  \{Z,Y,P  \}$ where  $Z \subseteq  Y 
\subseteq   P \subset H$,  and  $Z$   is   a   locally  complete    
intersection zero-dimensional  subscheme, and  $Y$,$P$ are  curves.  If we 
denote by $g$ the arithmetic  genus of $C$, by $d$,$y$,$p$ the degrees 
of $C$,$Y$,$P$ respectively, and by $z$ the length of $Z$, then 
\begin{align*} 
d&= y+p \\ 
g&= \frac{1}{2} (y-1)(y-2) + \frac{1}{2} (p-1)(p-2) +y-z-1. 
\end{align*} 
\ep 
 
The computation of the degree and genus of $C$ is straightforward from 
the exact sequence~(\ref{first}).

\br 
The curve $C$ is arithmetically Cohen-Macaulay if and only if $Z$ is  
empty. From the exact sequence~(\ref{first}) we obtain an exact sequence 
$$ 
0 \xrightarrow{} \mbox{H}^{1}_{*} (\ideal{C})  
\xrightarrow{} \mbox{H}^{1}_{*} ( \ideal{Z,H}(-p))  
\xrightarrow{} \mbox{H}^{2}_{*} (\ideal{Y} (-1) )
$$ 
where, for a a coherent sheaf $\caf$ on $\pso$, 
we let      
$$      
\mbox{H}^{i}_{*} ( \caf) = \bigoplus_{n \in \Z}       
\mbox{H}^{i} ( \caf (n) ).      
$$ 
If $Z$  is empty, the middle term  is zero, and that  forces the first 
term, which  is the Rao module $M_{C}$  of $C$ to be  zero also. Hence 
$C$ is ACM. 
 
Conversely, if $M_{C} = 0$, then we look at the particular twist 
$$ 
0 \xrightarrow{}  
\mbox{H}^{1}  (\ideal{Z,H} (-1)) \xrightarrow{}   
\mbox{H}^{2} (\ideal{Y} (p-2)). 
$$ 
The last group is isomorphic to  $\mbox{H}^{1} (\coo_{Y} (p-2))$, 
which is zero because $p \geq y$. Hence 
$\mbox{length} \,  Z = \dim \mbox{H}^{1}  (\ideal{Z,H} (-1))=0$,  
and $Z$ is  
empty. 
\er

\section{Existence of curves with a given triple} \label{three} 
 
Suppose  given $\{Z,Y,P\}$  as in  \ref{1.1}. Let  $\scrl$ be  the sheaf 
$\coo_{Y} (Z-1)$ associated  to the generalized divisor $Z$  on $Y$ as 
in \cite{hgd}, 2.8.  Then we will  show how to construct  a surjective 
map $u:  \ideal{P,\pso} \rrr \scrl$,  so as to define  $\ideal{C,\pso} = 
\mbox{Ker} \, u$ as in  the diagram~(\ref{second}). Note that $C$ will 
be locally Cohen-Macaulay and pure dimensional because at closed 
points $\scrl$ has depth one and $\ideal{P, \pso}$ has depth two.  
 
The  map $u$  needs  to  be compatible  with  the existing  inclusions 
$\coo_{\pso} (-1)  \hookrightarrow \ideal{P,\pso}$ and  $\coo_{Y} (-1) 
\hookrightarrow \scrl$. Furthermore,  $\scrl$ is an $\coo_{Y}$-module, 
so $u$ must factor  through $\ideal{P, \pso} \otimes \coo_{Y}$. Tensor 
the middle  line of the diagram~(\ref{second}) with  $\coo_{Y}$. Then it 
will be  sufficient to define  a surjective map $\bar{u}$  to complete 
the diagram 
\begin{equation} \label{third} 
\begin{CD} 
0 @>>>  \coo_{Y} (-1)  @>>> \ideal{P}/{\ideal{Y} \ideal{P}}   
@>>> \coo_{Y} (-p) @>>> 0  \\ 
   &&      @VVV            @VV\bar{u}V            @VVwV                \\ 
0 @>>> \coo_{Y} (-1) @>>>  \scrl       @>>>   \omega_{Z} (2-y)  @>>> 0 
\end{CD} 
\end{equation} 
Since  $P$  is  a  complete   intersection  in  $\pso$,  the  top 
row, which is the restriction to $Y$ of the conormal sequence of $P$, 
splits. Thus to  find $\bar{u}$, we need only find  a map $v: \coo_{Y} 
(-p) \rrr \scrl$ whose image in $ \omega_{Z} (2-y)$ is surjective. Since 
$Z$ is a locally complete  intersection in $H$, the sheaf $\omega_{Z}$ 
is generated by  a single element at every  point, so any sufficiently 
general   $w:   \coo_{Y}  (-p)   \rrr   \omega_{Z}   (2-y)$  will   be 
surjective.  To  show  that  $w$   lifts  to  $v$,  think  of  $v  \in 
\mbox{H}^{0} (Y, \scrl (p))$. Then we have an exact sequence 
$$ 
0 \xrightarrow{} \mbox{H}^{0} (\coo_{Y} (p-1))  
\xrightarrow{} \mbox{H}^{0} ( \scrl (p))  
\xrightarrow{} \mbox{H}^{0} (\omega_{Z} (2+p-y) ) 
\xrightarrow{} \mbox{H}^{1} (\coo_{Y} (p-1)) 
$$ 
Now  $p  \geq  y$,  so  $p-1  >  y-3$  hence  $\mbox{H}^{1}  (\coo_{Y} 
(p-1))=0$. Thus any surjective $w$  lifts to $v$ and gives the desired 
$\bar{u}$.  We  can  now  define  $C$  by  setting  $\ideal{C,\pso}  = 
\mbox{Ker} \,  u$. Since $\scrl \cong \ideal{P}/\ideal{C}$  has depth  
$\geq 1$ at  closed points, $C$ is locally Cohen-Macaulay,  and it is 
clear from our construction that $T(C) = \{Z,Y,P \}$. 
 
Finally,  we claim  that the  construction above  gives a  one  to one 
correspondence  between curves $C  \subset 2H$ 
with $T(C) = \{Z,Y,P \}$ and  global sections $ v \in \mbox{H}^{0} (Y, 
\scrl  (p))$  whose  image  in $\mbox{H}^{0}  (  \omega_{Z}  (2+p-y))$ 
generates $\omega_{Z} (2+p-y)$ at every  point. Note that it is enough 
to show that a surjective $\bar{u}$ fitting in the diagram~(\ref{third}) 
is determined by its  kernel.  Now, if $\bar{u}_{1}$ and $\bar{u}_{2}$ 
have the  same kernel and  fit in the diagram~(\ref{third}),  then their 
difference factors  through a  morphism $\omega_{Z} (2-y)  \rrr \scrl$ 
which must be zero, thus $\bar{u}_{1} = \bar{u}_{2}$. So we have: 
\bp \label{2.1} 
For each triple  $\{Z,Y,P \}$ as in~\ref{1.1} there  exists a curve $C 
\subset  2H$  with  $T(C) =  \{Z,Y,P  \}$.  The  set  of such  $C$  is 
parametrized  by an  open  subset of  the  vector space  ${\it H}^{0} 
(Y,\scrl (p))$, which has dimension 
$$ 
h^{0} ( \scrl (p)) = z + (p-1)y + 1 - \frac{1}{2} (y-1)(y-2). 
$$ 
\ep 
 
\section{Families of curves in the double plane} \label{four} 
 
In order to understand the Hilbert scheme of curves in the double 
plane, we must carry out the constructions of section \ref{two} and 
\ref{three} for flat families. Consider a base scheme $S$, and let $\ccc 
\subset 2H \times S$ be a family of curves of degree $d$ and genus 
$g$, i.e. a closed subscheme, flat over $S$, whose fibre over each 
point $s \in S$ is a locally Cohen-Macaulay curve $C_{s}$ in $2H$ of 
given degree and genus. The functor which to each $S$ associates the 
set of such flat families is represented by a Hilbert scheme which we 
denote by $H_{d,g}(2H)$. 
 
If $\ccc \subset 2H \times S$ is a flat family as above, the 
intersection $ \ccc \cap (H \times S)$ need not be flat. For example, 
if $S$ is integral, 
flatness of $ \ccc \cap (H \times S)$ is equivalent to local constancy 
of the integers $z$,$y$,$p$ associated to the fibres $C_{s}$ as in 
section~\ref{two} above. Applying Mumford's flattening stratification 
to $ \ccc \cap (H \times S)$, where $\ccc$ is the universal family 
over $H_{d,g}(2H)$, we find that the scheme $H_{d,g}(2H)$ is 
stratified by locally closed subschemes $H_{z,y,p} (2H)$ representing 
families for which $ \ccc \cap (H \times S)$ is flat and the fibres 
correspond to $z$,$y$,$p$ as above. 
 
We claim that the procedures of section~\ref{two} and \ref{three} 
above can be relativized over a base scheme $S$, to show that a flat 
family $\ccc \subset 2H \times S$ with the condition that $ \ccc \cap 
(H \times S)$ is also flat gives rise to a triple $\{Z,Y,P  \}$ where 
$Z \subseteq  Y \subseteq   P $ are closed subschemes of $H \times S$, 
flat over $S$, and where the fibres of $Z$ are zero-dimension locally 
complete intersection subschemes of $H$, and the fibres of  $Y$ and 
$P$  are  curves in $H$. Conversely, given any such triple  $\{Z,Y,P 
\}$, there are families of curves giving rise to this triple, and the 
resulting curves are parametrized by a global section of the sheaf 
$\scrl (P)$, where $\scrl$ is the sheaf  
$\Hom (\ideal{Z,Y} , \coo_{Y}(- 1))$ and $P$ is the divisor on $H 
\times S$ corresponding to $P$.  
 
Thus there is a natural map from the scheme $H_{z,y,p} (2H)$ to the 
Hilbert flag scheme $D_{z,y,p} (H)$ which parametrizes triples of flat 
families $Z \subseteq  Y \subseteq   P $ in $H \times S$, flat over 
$S$, as above, whose fibres $Z_{s}$ are locally complete intersection 
of length $z$, and $Y_{s}$ and $P_{s}$ are curves of degree $y$ and 
$p$ respectively. The fibres of $H_{z,y,p} (2H)$ over $D_{z,y,p} (H)$ 
are open subsets of the vector spaces  
$\mbox{H}^{0} (Y_{s}, \scrl_{s}  (p))$, which have constant dimension, and 
thus $H_{z,y,p} (2H)$ appears as an open subscheme of a geometric 
vector bundle over $D_{z,y,p} (H)$. 
 
The functorial details to justify all this are standard, if rather 
lengthy, so we leave them to the reader and content ourselves with 
summarizing the results in the following proposition.  
 
\bp \label{vect} 
The Hilbert scheme $H_{d,g} (2H)$ of curves in $2H$ of degree $d$ and 
genus $g$ is stratified by locally closed subschemes $H_{z,y,p} (2H)$ 
corresponding to families of curves whose integers $(z,y,p)$ 
associated by \ref{1.1} are constant. 
 
The scheme $H_{z,y,p}$ has a natural map to the Hilbert flag scheme 
$D_{z,y,p}(H)$ of triples $\{Z,Y,P  \}$ as in~\ref{1.1}, and this map makes 
$H_{z,y,p} (2H)$ into an open subset of a geometric vector bundle 
$\cae$ over $D_{z,y,p} (H)$, where $\cae$ is locally free of rank  
$$z + (p-1)y + 1 - \frac{1}{2} (y-1)(y-2).$$ 
\ep 
 
Next, we study the Hilbert flag scheme $D_{z,y,p} (H)$. 
\bp 
Given integers $z$,$y$,$p$ with $z \geq 0$, $p \geq y \geq 1$, the 
Hilbert flag scheme $\bar{D}_{z,y,p} (H)$ of closed subschemes  
$Z \subseteq  Y \subseteq   P $ in $H$, with $Z$ of dimension zero and 
length $z$, and $Y$, $P$ curves of degree $y$,$p$ respectively, is 
irreducible and generically smooth of dimension 
$$ 
z + \frac{1}{2} \, y (y+3) + \frac{1}{2} \, (p-y) (p-y+3). 
$$ 
\ep 
\begin{proof} 
Since $Y \subseteq P$, we can write $P=Y+W$, where $W$ is an effective 
divisor,  i.e. a  curve, of  degree  $p-y$, which can be chosen 
independently of $Z$,$Y$. So 
$$\bar{D}_{z,y,p}(H) = D_{z,y} (H)  \times D_{p-y}(H), $$ where $D_{p-y}$ is 
a   projective   space   of   dimension   $\frac{1}{2}   (p-y)(p-y+3)$ 
parametrizing curves  of degree $p-y$.  Thus we reduce  to considering 
$D_{z,y} (H)$. This is the  non-trivial part, which was proved by Brun 
and Hirschowitz \cite{brun-hirsch}, proposition 3.2. The main ideas of 
their proof are as follows.   
 
We regard $D_{z,y}$ as a closed subscheme 
of  $M \times  K$,  where $M=  Hilb^{z}  (H)$ is  the  scheme of  zero 
dimensional  closed subschemes  of $H$  of length  $z$, and  $K$  is a 
projective  space  of dimension  $\frac{1}{2}  y (y+3)$  parametrizing 
curves of  degree $y$ in $H$.  Let $Z$ and $Y$ denote the universal
families over $M$ and $K$ respectively. It is a theorem  of Fogarty \cite{hl}, 
example 4.5.10,  
that $M$ is smooth irreducible of dimension $2z$. 
Brun and Hirschowitz first show that 
$D_{z,y}$  is the  zero scheme  of a  section of  the rank  $z$ vector 
bundle $\Hom  (\ideal{Y} (y), \coo_{Z} (y))$  on $M \times  K$, and so 
has  codimension at most  $z$ at  every point.  Thus the  dimension of 
$D_{z,y}$ is at least $z +  \frac{1}{2} y (y+3)$ at every point. Next, 
observe that  $D_{z,y}$ contains an  open subset $U$  corresponding to 
divisors of degree $z$ on smooth curves, that is smooth irreducible of 
dimension $z +\frac{1}{2} y(y+3)$. To complete the proof, they use the 
theorem  of  Brian{\c{c}}on  \cite{brian}  that the  punctual  Hilbert 
scheme  of zero  dimensional  schemes of  length  $z$ at  a point  has 
dimension $z-1$. This shows that the fibre of $D_{z,y}$ over any curve 
$Y$ has  total dimension at most  $z$. Thus the  dimension of $D_{z,y} 
\setminus U$ must be strictly  less than $z+ \frac{1}{2} y (y+3)$, and 
hence it is contained in the closure of $U$. 
\end{proof} 
 
\bc  \label{dim} 
Suppose  $z \geq 0$ and $p \geq y \geq 1$, or $z=y=0$ and $p \geq 1$. 
Then the scheme $H_{z,y,p} (2H)$ is 
irreducible and generically smooth of dimension  
$$ 2z + \frac{1}{2} y (y+1)+\frac{1}{2} p (p+3).$$ 
\ec 
 
\section{Irreducible components of $H_{d,g} (2H)$} \label{fours} 
 
We can now describe the irreducible components of $H_{d,g} (2H)$. 
\bt \label{components} 
Let $d$ and $g$ be integers, with $d \geq 1$. 
$H_{d,g} (2H)$ is non-empty if and only if either   
$g = \frac{1}{2} (d-1)(d-2)$, or $d \geq 2$ and  
$g \leq \frac{1}{2} (d-2)(d-3)$.   
 
The irreducible components of $H_{d,g}$ are the closures 
$\bar{H}_{z,y,p}$ of the subschemes $H_{z,y,p}$ defined in 
section~\ref{four}, where $(z,y,p)$ varies in the set  
of triples of nonnegative integers  satisfying the following conditions: 
$p \geq 1$, $p \geq y$ ,$z=0$ if $y=0$, and 
\begin{align*} 
p &= d-y \\ 
z &= \frac{1}{2} (y-1)(y-2) + \frac{1}{2} (p-1)(p-2) +y-g-1. 
\end{align*} 
\et 
 
\begin{proof} 
By corollary \ref{dim}, $\bar{H}_{z,y,p}$ is 
irreducible of dimension $2z + \frac{1}{2} y (y+1)+\frac{1}{2} p (p+3)$.  
To see 
that $\bar{H}_{z,y,p}$ is an irreducible component, it is enough to 
observe that its generic point is not contained in a different 
$\bar{H}_{z',y',p'}$. Now this is clear, since the support of the 
generic  point of $\bar{H}_{z,y,p}$ is the union of two  smooth plane  
curves of degrees $y$ and $p-y$ respectively.  
 
By propositions~\ref{1.1} and \ref{2.1}, as a topological space 
$H_{d,g}$  is the disjoint union of its subschemes $H_{z,y,p}$, where 
$z$,$y$,$p$ are nonnegative integers satisfying  
$p \geq 1$, $p \geq y$, $z=0$ if $y=0$, and 
\begin{align*} 
p &= d-y \\ 
z &= \frac{1}{2} (y-1)(y-2) + \frac{1}{2} (p-1)(p-2) +y-g-1. 
\end{align*} 
It remains to determine the triples $(z,y,p)$ satisfying these 
conditions. We must have   
$$z=\frac{1}{2} (d-2)(d-3)-g - (y-1)(d-y-2).$$ 
We now impose the conditions on $z$, $y$ and $p$. The conditions  
$p \geq y$ and $p \geq 1$ translate into $y \leq d/2$ and $d-y -1 
\geq 0$. As $y$ and $z$ must be nonnegative, we have 
$$\frac{1}{2} (d-1)(d-2) -g = y(d-y-1)+z \geq 0$$ 
with equality if and only if $y=z=0$ or $z=0$, $y=1$ and 
$d=2$. In all other cases we have $y \geq 1$ and either $d-y >1$, or 
$y=1$,$d=2$ and $z \geq 1$, hence  
$$\frac{1}{2} (d-2)(d-3) -g = (y-1)(d-y-2)+z \geq 0.$$  
\end{proof} 
 
\br \label{rcomponents}  
How many irreducible components does $H_{d,g}$ have ?  
If $d \neq 2$ and $g=1/2(d-1)(d-2)$, $H_{d,g}=H_{0,0,d}$ is 
irreducible. $H_{2,0}$ has two irreducible components, namely 
$\bar{H}_{0,1,1}$ and $\bar{H}_{0,0,2}$. 
 
Suppose $d \geq 2$ and $g \leq 1/2 (d-2)(d-3)$ with $(d,g) \neq 
(2,0)$. Let $y_{M}$ be largest integer $n$ in 
the closed interval $[1,d/2]$ such that  $(n-1)(d-n-2) \leq  
1/2 (d-2) (d-3) - g$. The irreducible components of $H_{d,g}$  
are in one to one correspondence with the integers $y$ in the closed  
interval $[1,y_{M}]$. 
\er

\section{Numerical invariants}  \label{five} 
Out of the exact sequence $(\ref{first})$ one immediately computes the 
postulation of the curve $C$ in  terms of the two integers $y$ and $p$ 
and of  the postulation of $Z$.  However to compute  the dimensions of 
the first and  second cohomology modules of $\ideal{C} (n)$ one needs 
some additional  argument. For example,  one could use the  results on 
liaison of  the next  section. Here we  take a different  approach. 
 
Choose a point $R$ of $\pso \setminus H$, and let $\pi: 2H \rrr H$ be 
the morphism induced by the projection from $R$. $\pi$ is a finite and 
flat morphism, and $\pi_{*} (\coo_{2H}) = \coo_{H} \oplus \coo_{H} (-1)$. 
Now given any other point $R'$ in $\pso \setminus H$, there is an 
automorphism of $\pso$ which sends $R$ to $R'$ and fixes every point 
of $H$. It follows that for any coherent sheaf of $\coo_{2H}$-modules 
$\caf$, the isomorphism class of the sheaf  $\pi_{*} \caf$ is 
independent of the choice of the point $R$.

\bp \label{cae} 
Let $\pi: 2H \rrr H$ be 
the morphism induced by the projection from $R$. 
Let $C$ be a curve in $2H$. Then  
$\cae = \pi_{*} (\ideal{C,2H})$ is a locally free $\coo_{H}$-module 
of rank two. Furthermore, the Rao module  
\begin{em} 
$\mbox{H}^{1}_{*} (\pso, \ideal{C})$ 
\end{em} 
is  isomorphic to  
\begin{em} 
$\mbox{H}^{1}_{*} (H, \cae)$  
\end{em} as a module over the  
homogeneous coordinate ring of $\pso$.  
\ep 
\begin{proof} 
We use the fact proven in \ref{1.1} that the sheaf of 
$\coo_{Y}$-modules $\scrl$, 
defined as the quotient $\ideal{P} / \ideal{C}$ has depth $\geq 
1$ at closed points. It follows from this and the exact sequence 
$$ 
0 \rrr \pi_{*} (\ideal{C,2H}) \rrr 
\pi_{*} (\ideal{P,2H}) \rrr \scrl \rrr 0 
$$ 
that $\pi_{*} (\ideal{C,2H})$ is locally free of rank two provided 
$\pi_{*} (\ideal{P,2H})$ is locally free of rank two . But this is 
clear as we have an exact sequence 
$$ 
0 \rrr \coo_{H}(-1) \rrr 
\pi_{*} (\ideal{P,2H}) \rrr \ideal{P,H} \cong \coo_{H} (-p) \rrr 0. 
$$  
 
Finally, from the exact sequence $(\ref{first})$ 
we see that the Rao module of $C$ is annihilated by the equation of 
$H$, hence it is isomorphic to  
$\mbox{H}^{1}_{*} (H, \pi_{*} (\ideal{C,2H}) )$ as a modules over 
$\mbox{H}^{0}_{*} (\pso, \coo_{\pso} )$. 
\end{proof}

\bc 
Let $C$ be a curve of degree $d$ in the double plane, with Rao module 
$M_{C}$. Let $M_{C}^{*}$ denote the $k$-dual module to $M_{C}$ (see 
\cite{MDP}, 0.1.7 p.20). Then $M_{C}^{*} \cong M_{C} (d-2)$, and in 
particular $h^{1} (\ideal{C} (n))=h^{1} (\ideal{C} (d-2-n)$ for all 
integers $n$. 
\ec 
\begin{proof} 
Let $M=M_{C} \cong \mbox{H}^{1}_{*} (H, \cae)$. Since $\cae$ is a rank 
two vector  bundle on $H$,  we have $\dual{\cae} \isom  \cae (-c_{1})$ 
where  $c_{1}$ denotes  the  first  Chern class  of  $\cae$. By  Serre 
duality  we  have  $M^{*}  \cong  M(-c_{1}-3)$.  Now  from  the  exact 
sequence 
\begin{equation} \label{caeseq} 
\exact{\coo_{H} (-y-1)}{\cae}{\ideal{Z,H}(-p)}. 
\end{equation} 
we compute $c_{1} = -p-y-1=-d-1$, and we are done.  
\end{proof} 
 
\br 
If we apply $\pi_{*}$ to the exact sequence $(\ref{first})$, we obtain 
an extension class in $\mbox{Ext}^{1}_{H} (\ideal{Z,H} 
(-p),\coo_{H}(-y-1))$, hence an element $w$ in $\mbox{H}^{0}  
(\omega_{Z} (2+p-y))$. In the notation of section~\ref{three}, one 
verifies that $w$ is the image of the global section $v$ of $\scrl 
(p)$  corresponding to $C$, under the map  
$\mbox{H}^{0} (Y, \scrl (p)) \rrr \mbox{H}^{0}  
(Y, \omega_{Z} (2+p-y))$ coming from diagram~\ref{third}. 
\er 
\br  
We may ask which locally free sheaves $\cae$ arise from this 
construction. Let us say that a locally free sheaf $\caf$ is normalized if 
it has a global section, while $\caf (-1)$ has no section. Let $\caf$ 
be a locally free normalized $\coo_{H}$-module of rank two. Then one 
may check that there exists a curve $C \subset 2H$ such that 
$\pi_{*} (\ideal{C,2H})$ is 
isomorphic to $\caf(n)$ for some integer $n$, if and only if $c_{1} 
\caf \leq 1$. In particular, all unstable locally free sheaves of rank 
two arise from this construction. 
\er

\bc \label{cohomology} 
Let $C$ be a curve in $2H$ with $T(C)= \{Z,Y,P\}$. Then 
\begin{align*} 
h^{0} ( \ideal{C} (n)) &= 
h^{0} ( \coo_{\pso} (n-2)) +  
h^{0} ( \coo_{H} (n-y-1))  +\\ 
& + h^{0} ( \ideal{Z,H} (n-p)) 
\end{align*} 
and 
\begin{align*} 
h^{2} (\ideal{C} (n)) & = 
h^{0} (\coo_{\pso} (-n-4)) - h^{0} ( \coo_{\pso} (-n-2)) + \\ 
& + h^{0} (\coo_{H} (p-3-n)) + 
h^{0} (\ideal{Z,H} (y-2-n)). 
\end{align*} 
\ec 
\begin{proof} 
The formula for the postulation follows immediately from the exact 
sequence $(\ref{first})$. To obtain the formula for $h^{2}$ we  
compute $h^{2} \cae (n)$ using Serre duality:  
$$ 
h^{2} \cae (n) = h^{0} \cae (p+y-2-n) = h^{0} \coo_{H} (p-3-n) + 
h^{0} \ideal{Z,H} (y-2-n). 
$$ 
\end{proof} 
 
We can also determine the postulation character of a curve $C$ in 
$2H$. Recall \cite{MDP} for any numerical function $f(n)$ one defines 
the difference function $\df f(n) = f(n) - f(n-1)$. If $C$ is a curve 
in $\pso$, one defines its {\em postulation character} $\gmc$ by 
$$ 
\gmc (n) = \df^{3} ( h^{0} ( \ideal{C} (n)) -  h^{0} ( \coo_{\pso} (n)) ). 
$$ 
Similarly, for a zero-dimensional closed subscheme $Z \subset \ptwo$ 
we define its postulation character 
$$ 
\gmz (n) = \df^{2} ( h^{0} ( \ideal{Z} (n)) -  h^{0} ( \coo_{\ptwo} (n)) ). 
$$ 
By analogy with the case of ACM curves \cite{MDP}, I.2 and V.1.3, one 
shows easily the following: 
\bp \label{gamma} 
Let $Z$ be a zero dimensional closed subscheme of $\ptwo$ of degree 
$z$, and let $s$ be the least degree of a curve containing $Z$. Then 
$$  
\gmz (n) 
= \left\{  
\begin{array}{ll} 
0 &   \mbox{ if \;\;\; $ n < 0$},       
\\      
-1 
&       
\mbox{ if \;\;\; $ 0 \leq n < s$},       
\\             
a_{n} \geq 0      
&       
\mbox{ if \;\;\; $ n \geq s$}. 
\end{array}       
\right.       
$$ 
Furthermore, $\gmz$ is a character, so $\sum_{n \geq s} a_{n} = s$, 
and we can determine the degree z by 
$$ 
z= \sum_{n \geq s} n a_{n} - \frac{1}{2} s (s-1). 
$$ 
 
Conversely, given integers $s \geq 1$ and $a_{n} \geq 0$ for $n \geq 
s$ there exists a reduced zero-dimensional closed subscheme $Z \subset \ptwo$ 
with postulation character as above.  
\ep 
 
Now using (\ref{cohomology}) we can compute the postulation character 
of $C$. We find 
\bt \label{gammac} 
Let $C$ be a curve in $2H$ with associated triple $\{Z,Y,P \}$. Then 
there is an integer $s$ with $0 \leq s \leq y$ and there are integers 
$b_{n}$ for $n \geq p+s$ such that $\sum b_{n} = s$, and the 
postulation character $\gmc$ is the sum of the following functions: 
$-1$ in degree $0$, $-1$ in degree $1$, $1$ in degree $y+1$, $1$ in 
degree $p+s$, and $\df \beta$, where 
$$ 
\beta (n) = \left\{  
\begin{array}{ll} 
0 &   \mbox{ if \;\;\; $ n < p+s$}       
\\      
b_{n} 
&       
\mbox{ if \;\;\; $ n \geq p+s$}. 
\end{array}       
\right.  
$$ 
 
Conversely, given integers $0 \leq s \leq y \leq p$ and $b_{n} \geq 0$ 
for $n \geq p+s$ with $\sum b_{n} =s$, there exists a curve $C$ in 
$2H$ with postulation character as described. 
\et 
 
\br 
In the notation of (\ref{gammac}), the degree $z$ of $Z$ is determined 
by  
$$ 
z = \sum_{n \geq p+s} (n-p) b_{n} - \frac{1}{2} s (s-1). 
$$ 
\er 
 
\br 
One can see easily that the possible postulation characters of curves 
on surfaces of degree $2$ other than $2H$ form a subset of the 
postulation characters described in (\ref{gammac}). Thus we have found 
all postulation characters of curves in surfaces of degree $2$. This 
gives some hope that one day it will be possible to determine all 
possible postulation characters of curves in $\pso$. 
\er

\section{Liaison} \label{six} 
In this section we use notation and terminology of \cite{hgd}, section 
4.  
The following proposition describes the behaviour of the triple  
$T(C) = \{Z,Y,P\}$ under liaison. 
\bp 
Let $C$ be a curve in $2H$ with $T(C)= \{Z,Y,P\}$.  
Suppose $S$ is a surface containing $C$ and meeting $H$ properly.  
Let $Q= H \cap S$. Then the curve $D$ linked to $C$ by the complete 
intersection $2H \cap S$ has triple $T(D)=\{Z,Q-P,Q-Y\}$; in particular, 
$Z \subset Q-P$. 
 
Conversely, if $Q \subset H$ is a curve containing $P$ such that  
$Z \subset Q-P$, then $C$ is linked by the complete intersection of 
$2H$ with some surface $S$ to a curve $D$ with $T(D)=\{Z,Q-P,Q-Y\}$. 
 
In particular,  $Z$ and $P-Y$ are 
invariant under liaison on $2H$. 
\ep 
\begin{proof} 
Suppose first $C$ is linked to $D$ by the complete intersection $E= 2H 
\cap S$, and let $T(D)=\{Z',Y',P' \}$. Let $Q$ be the intersection of 
$S$ with the reduced plane $H$. We claim that $Z'=Z$, $Y'=Q-P$ and 
$P'=Q-Y$.  
 
By \cite{hgd} page 317 the ideal sheaf of $D$ in $2H$ is 
$\Hom_{\coo_{2H}} ( \ideal{C,2H}, \ideal{E,2H})$ with its natural 
embedding in $\Hom_{\coo_{2H}} ( \ideal{E,2H}, \ideal{E,2H}) \cong 
\coo_{2H}$. By \cite{AG} exercises III 6.10 and 7.2, applying the functor 
$\Hom_{\coo_{2H}} ( -, \ideal{E,2H})$ to the exact sequence 
$$ 
 0 \rrr  
\ideal{Y,H}(-1) \stackrel{h}{\rrr} \ideal{C,2H} \rrr  
\ideal{Z,H} (-P)  \rrr 0 
$$ 
we obtain a long exact sequence 
$$ 
0 \rrr  
\ideal{Q-P,H}(-1) \stackrel{h}{\rrr} \ideal{D,2H} \rrr 
\ideal{Q-Y,H} \rrr \omega_{Z} (p-q+2) \rrr 0. 
$$ 
As $\omega_{Z} (p-q+2) \cong \coo_{Z} (y-q)$, the kernel of the last 
map is $ \ideal{Z,H}(Y-Q)$, and our claim follows. 
 
Conversely, suppose that $Q$ is a curve in $H$ containing $P$, and 
such that $Z \subset Q-P$. The exact sequence~(\ref{first}) tells us that 
there exists a surface $S$ containing $C$ but not $H$, whose 
intersection  with $H$ is $Q$. By what we have just shown,  
$T(D)=\{Z,Q-P,Q-Y\}$. 
\end{proof} 
 
\bc 
Let $C$ be a curve in $2H$ with $T(C)= \{Z,Y,P\}$, and let $Y'$ be a 
curve in $H$ containing $Z$. Let $y$ and $y'$ be the degrees of $Y$ 
and $Y'$ respectively. There is curve $D$ obtained from $C$ 
by an elementary biliaison of height $y'-y$ on $2H$ with $T(D)= 
\{Z, Y',Y'+P-Y \}$.   
\ec 
\begin{proof} 
Use the above proposition with $Q=P+Y'=(Y'+P-Y) + Y$. 
\end{proof} 
 
\bc 
Let $C$ be a curve in $2H$ with $T(C)= \{Z,Y,P\}$. Suppose that $C$ is 
not arithmetically Cohen-Macaulay, that is, $Z$ is not empty. Then $C$ 
is minimal in its biliaison class if and only if $Y$ has minimal 
degree among curves in $H$ containing $Z$. 
\ec 
\begin{proof} 
If $h^{0} (H, \ideal{Z,H} (\deg Y -1) > 0$, by the previous corollary 
there is a curve obtained from $C$ by an elementary biliaison of 
negative height, hence $C$ is not minimal. 
 
If $h^{0} (H, \ideal{Z,H} (\deg Y -1) = 0$, then by 
corollary~\ref{cohomology} we have  
$$h^{2} (\pso, \ideal{C} (1)) - 2 h^{2}  (\pso, \ideal{C} ) +  
h^{2} (\pso, \ideal{C} (-1) ) \leq 1.$$ Now it follows from \cite{MDP} 
Proposition~III.3.5 that $C$ is minimal: see \cite{sch-tran} Corollary~4.4. 
\end{proof} 
 
\section{Connectedness of the Hilbert scheme} \label{seven} 
\bt \label{connect} 
The Hilbert scheme $H_{d,g} (2H)$ is connected. 
\et 
\begin{proof} 
If $g=\frac{1}{2} (d-1)(d-2)$ and $d\neq 2$, $H_{d,g}$ is irreducible by  
theorem~\ref{components}.  
 
To handle the case  $d=2$ and $g=0$, we fix homogeneous  
coordinates $[x:y:z:w]$ on $\pso$, so that $x=0$ is an equation for 
$H$, and we look at the family of curves $C_{t}$ in $\pso \times Spec \, 
k[t]$ defined by the global ideal 
$$ 
I = < x^{2}, xy, y^{2}, x+ty >. 
$$ 
For $t\neq 0$, $C_{t}$ is in $H_{0,1,1}$, while $C_{0}$ belongs to 
$H_{0,0,2}$. It now follows from remark~\ref{rcomponents} that 
$H_{2,0}$  is connected. 
 
If $g < \frac{1}{2}(d-1)(d-2)$, by remark~\ref{rcomponents} we have set 
theoretically  
$$ 
H_{d,g} = \bigcup_{1 \leq y \leq y_{M}} H_{z_{y},y,d-y} 
$$ 
where  
$$ 
z_{y} = \frac{1}{2} (d-2) (d-3) - g - (y-1)(d-y-2). 
$$ 
In particular, $H_{d,g}$ is irreducible for $d \leq 3$, and for    
$d \geq 4$ the theorem is a consequence of the following 
proposition. 
\end{proof} 
\bp \label{pconnect} 
If $y \geq 2$, $p \geq y$ and $r \geq 0$,  
there is a curve in $H_{r,y,p}$ 
specializing to one in $H_{r+p-y,y-1,p+1}$. 
\ep 
\begin{proof} 
Suppose the claim is true when $y=2$. Then 
by adding $y-2$ times a plane section (i.e. by performimg an 
elementary biliaison of height $y-2$ on $2H$, see \cite{hmdp1} 
proposition~1.6) we find that the claim is true for all $y \geq 2$, 
and we are done. 
 
We now construct a family of curves in $H_{r,2,p}$ 
specializing to one in $H_{r+p-2,1,p+1}$. This has already been done 
by Nollet in \cite{nthree}, example 3.10, in the case $p=2$ and $r=1$, 
and his construction,  generalizes without any major 
modification (Nollet noticed this independently \cite{nolletp}). 
Here are the details.  
 
As above we fix homogeneous coordinates $[x:y:z:w]$ on $\pso$, so that 
$x=0$ is 
an equation for $H$. We let $Y$ and $P$ denote respectively the double 
line  $x=y^{2}=0$ and the multiple line $x=y^{p}=0$. Let $Z$ have equations 
$x=y=f=0$ where $f$ is a form of degree $r = \deg Z$ in $k[z,w]$.  
 
To give a curve $C$ in $2H$ with $T(C) = \{Z,Y,P \}$ is by 
proposition~\ref{2.1} the same as giving a morphism $v : 
\ideal{Z,H} (-p) \rrr \coo_{Y} (-1)$ whose image in  in $\mbox{H}^{0}  
(\omega_{Z} (p))$ generates $\omega_{Z} (p)$ at every point. This 
amounts to choosing forms $s$ and $g$ in $k[z,w]$, of degrees $p-1$ and 
$r+p-2$ respectively, and such that $g$ and $f$ have no common zeros 
on the line $x=y=0$. The corresponding morphism $\ideal{Z,H} (-p) \rrr 
\coo_{Y} (-1)$ sends $y$ to $ys$ and $f$ to $fs+yg$. For the corresponding 
curve $C$ we have: 
$$ 
I_{C} = \,\, < x^{2}, xy^{2}, y^{p+1}+xys, xfs+xyg+y^{p}f >. 
$$ 
 
Now we look at the family of curves in $\pso \times Spec \, k[t]$ 
obtained by flattening the ideal generated by $x^{2}, A, B ,C$ where 
\begin{align*} 
A & =  xy^{2} \\ 
B & =  ty^{p+1}-xyz^{p-1} \\ 
C & =  xyw^{r+p-2}-tz^{r} (ty^{p}-xz^{p-1}). \\ 
\end{align*} 
For $t \neq 0$ we obtain a curve $C$ as above with $s=-\frac{1}{t} z^{p-1}$, 
$g=w^{r+p-2}$ and $f =-t^{2} z^{r}$. To see what happens at $t=0$, we 
set 
\begin{align*} 
D & =  \frac{1}{t} (w^{r+p-2} B+z^{p-1}C )=y^{p+1}w^{r+p-2}+xz^{r+2p-2}+tF \\ 
E & =  \frac{1}{t} (z^{p-1} A + y B) = y^{p+2}. \\ 
\end{align*} 
 
It follows that the ideal of the limit scheme $C_{0}$ contains the 
ideal  
$$ 
J= \,\, < x^{2}, xy^{2}, xyz^{p-1}, xyw^{r+p-2},y^{p+2}, 
y^{p+1}w^{r+p-2}+xz^{r+2p-2} >. 
$$ 
The saturation of $J$ is the ideal 
$$ 
I = \,\,< x^{2}, xy ,y^{p+2}, y^{p+1}w^{r+p-2}+xz^{r+2p-2} >. 
$$ 
But this is the homogeneous ideal of a curve $D$ in the double plane:  
$Y(D)$ is the line $x=y=0$, $P(D)$ has equations $x=y^{p+1}=0$, and 
$Z(D)$ is defined on $Y(D)$ by the equation $w^{r+p-2}=0$. Hence $D$ 
belongs to $H_{r+p-2,1,p+1}$. In particular, $D$ has the same degree 
and genus as $C_{0}$. So we must have $D=C_{0}$, and this finishes the 
proof.  
 
We remark that in this family the zero dimensional scheme $Z$ 
associated to $C_{t}$ is supported at the point $[0:0:1:0]$ for 
$t \neq 0$, and at the point $[0:0:0:1]$ for $t=0$ !  
\end{proof}

\section{Extremal and subextremal curves} 
 
Given a curve $C$, the function $\rho_{C}: \Z \rrr \Z$ defined by  
$\rho_{C} (n) =h^{1} \ideal{C} (n)$ 
is called the Rao function of $C$. It is the Hilbert function of the 
Rao module $M_{C} = \mbox{H}^{1}_{*} (\pso, \ideal{C})$ 
of $C$.

\bt[\cite{bounds},\cite{subextremal}] \label{bounds} 
Let $C \subset \pso$ be a curve of degree $d$ and arithmetic  
genus $g$. Then 
\begin{enumerate} 
\item 
$C$ is planar if and only if $g = \frac{1}{2} (d-1)(d-2)$.  
\item  
If $C$ is not planar, then $d \geq 2$, $g \leq \frac{1}{2} (d-2)(d-3)$ 
and  
$$ 
\rho_{C} (n) \leq \rho^{E}_{d,g} (n) \;\;\; \mbox{for all $n \in \Z$} 
$$   
where  
$$ 
\rho^{E}_{d,g} (n) =  
\left\{  
\begin{array}{ll} 
\mbox{\em max} (0, \rho^{E}_{d,g} (n+1)-1) &    
\mbox{ if \;\;\; $ n \leq -1$},       
\\      
\frac{1}{2} (d-2)(d-3) - g 
&       
\mbox{ if \;\;\; $ 0 \leq n \leq d-2$},       
\\  
\mbox{\em max} (0, \rho^{E}_{d,g} (n-1)-1)      
&       
\mbox{ if \;\;\; $ n \geq d-1$}. 
\end{array}       
\right.  
$$ 
\end{enumerate} 
\et 
 
\bd 
A curve $C \subset \pso$ of degree $d$ and genus $g$ is called {\em 
extremal}  if 
it is not planar and $\rho_{C} = \rho^{E}_{d,g}$. 
\ed 
\br 
The definition of an  
extremal curve in \cite{extremal},  
\cite{subextremal}  is different from ours, since  
they require the curve not to be ACM, but we allow the  
ACM case if $\rho^E_{d,g} =0$. 
\er 
\br 
There exist extremal curves (for example in $2H$, see below) if $d=2$ 
and $g \leq -1$ and if $d \geq 3$ and $g \leq \frac{1}{2} (d-2)(d-3)$.  
\er 
 
The following characterization of extremal curves is essentially 
contained in \cite{extremal} and \cite{subextremal}. 
 
\bp \label{extremal} 
Let $C$ be a curve of degree $d \geq 2$ and genus $g$.  
The following are equivalent: 
\begin{enumerate} 
\item 
$C$ is extremal; 
\item 
either $C$ is a minimal curve contained in the union of two distinct 
planes, or $C$ is obtained from a plane curve by an elementary  
biliaison of height one on a quadric surface; 
\item 
either $C$ is  ACM with $(d,g) \in \{ (3,0),(4,1) \} $ and $C$ is 
contained in an integral quadric surface, or $C$ is a non-planar 
curve which contains a plane curve of degree $d-1$. 
\end{enumerate} 
\ep 
\begin{proof} 
Looking at  \cite{bounds}, \cite{subextremal} one sees that $C$ is 
extremal if and only if it is not planar and $h^{0} \ideal{C} (2) \geq 
2$. From this we deduce that 1 implies 3. Curves $C$ in the union of two 
distinct planes $H_1 \cup H_2$ are studied  in \cite{hgd}, section  5: 
if $C \cap H_{1}$ has degree $d-1$, then $C$ is either minimal or 
ACM. From this we see that 3 implies 2.  
 
Finally, we prove 2 implies 1. If $C$ is obtained from a plane curve  
by an elementary  biliaison of height one on a quadric surface, then 
one computes $g = \frac{1}{2} (d-2)(d-3)$ and so $C$ is extremal. 
If $C$ is a minimal curve in  $H_1 \cup H_2$, then by \cite{hgd}, 
section  5, adding a suitable number  
(namely, $\frac{1}{2}(d-2)(d-3)-g-1 $)  of plane sections to $C$ we 
obtain a curve linearly equivalent to the 
disjoint union of two plane curves. Thus the Rao module of $C$ is 
isomorphic to $R/(h,k,f,g) \, (a-1)$ where $a=\frac{1}{2} (d-2)(d-3)-g$,  
$R$ is the homogeneous coordinate ring of $\pso$, $h$, $k$ are the 
equations of $H_{1}$ and $H_{2}$, and  $f$,$g$ are forms of degrees $a$, 
$d+a-2$ respectively, having no common zeros along the line 
$h=k=0$. Hence $\rho_{C} = \rho^{E}_{d,g}$. 
\end{proof} 
 
\bt[\cite{subextremal}] \label{subbounds} 
Let $C \subset \pso$ be a curve of degree $d$ and genus $g$, which is 
neither planar nor extremal. Then $d \geq 3$ and 
\begin{enumerate} 
\item 
$g \leq \frac{1}{2} (d-3)(d-4) + 1$; if equality holds, then $d \geq 
5$ and $C$ is ACM; 
\item 
if $d \geq 4$, then 
$$ 
\rho_{C} (n) \leq \rho^{S}_{d,g} (n) \;\;\; \mbox{for all $n \in \Z$} 
$$   
where  
$$ 
\rho^{S}_{d,g} (n) =  
\left\{  
\begin{array}{ll} 
\mbox{\em max} (0, \rho^{S}_{d,g} (n+1)-1) &    
\mbox{ if \;\;\; $ n \leq 0$},       
\\      
\frac{1}{2} (d-3)(d-4) + 1 - g 
&       
\mbox{ if \;\;\; $ 1 \leq n \leq d-3$},       
\\  
\mbox{\em max} (0, \rho^{S}_{d,g} (n-1)-1)      
&       
\mbox{ if \;\;\; $ n \geq d-2$}. 
\end{array}       
\right.  
$$ 
\end{enumerate} 
\et 
 
\bd 
A curve $C \subset \pso$ of degree $d \geq 4$ and genus $g$ is called 
{\em subextremal} if it is neither planar nor extremal and $\rho_{C} = 
\rho^{S}_{d,g}$. Again, this differs from the terminology of 
\cite{subextremal} in that we include among subetxremal curves  those 
ACM curves of degree $d \geq 5$ which have genus $1/2 (d-3)(d-4) +1$. 
\ed 
 
\br 
There exist subextremal curves (for example in $2H$, see below) if $d=4$ 
and $g \leq 0$ and  
if $d \geq 5$ and $g \leq \frac{1}{2} (d-3)(d-4)+1$.  
\er

\bp \label{subextremal} 
Let $C \subset \pso$ be a curve of degree $d \geq 4$ and genus $g$.  
The following are equivalent: 
\begin{enumerate} 
\item 
$C$ is subextremal; 
\item 
$C$ is obtained from an extremal curve by an elementary biliaison of 
height one on a quadric surface; 
\item 
either $C$ is ACM and $(d,g) \in  \{(5,2),(6,4) \}$, or $C$ is a 
divisor of type $(1,3)$ on a smooth quadric surface, or there is a plane $H$ 
such that $\ideal{C \cap H,H} = \ideal{Z,H} (2-d)$ with $Z$ zero-dimensional 
contained in a line, and the residual scheme $Res_{H} (C)$ to the 
intersection of $C$ with $H$ is a plane curve of degree two. 
\end{enumerate} 
\ep 
\begin{proof} 
If $C$ is not ACM, the equivalence of $1$ and $2$ is the content of 
Theorem~2.14 in \cite{subextremal}. The same proof works for the ACM case 
as well (cf. \cite{subextremal}, Lemma~2.5). 
 
We claim that $2$ implies $3$: suppose $C$ is obtained from an 
extremal  curve $B$ by an elementary biliaison of height one on a 
quadric surface $F$. 
If there is a plane $H$ whose intersection with $B$ has degree $d-3 =  
\deg  B -1 $ and if $F$ contains $H$ as a component, then   
$\ideal{C \cap H,H} = \ideal{Z,H} (2-d)$ with $Z$ zero-dimensional 
contained in a line, and $Res_{H} (C)$ is a plane curve. 
Otherwise, by proposition~\ref{extremal} either $B$ is an ACM curve with 
$(d,g) \in  \{(3,0),(4,1) \}$, or $B$ is a divisor of type $(0,2)$ on 
the smooth quadric surface $F$. So  $C$ is either ACM with  
$(d,g) \in  \{(5,2),(6,4) \}$, or $C$ is a divisor of type $(1,3)$ on 
$F$. 
 
Finally, we show that $3$ implies $1$. The special cases are clear, so 
we may assume that  there is a plane $H$ 
such that $\ideal{C \cap H,H} = \ideal{Z,H} (2-d)$ with $Z$ zero-dimensional 
contained in a line, and $Y = Res_{H} (C)$ is a plane curve of degree two. 
Then there is an exact sequence: 
\begin{equation} \label{last} 
 0 \rrr  
\ideal{Y,\pso}(-1) \stackrel{h}{\rrr} \ideal{C,\pso} \rrr  
\ideal{Z,H}(2-d)  \rrr 0  
\end{equation} 
from which we deduce that $ \mbox{length} \, Z 
=\frac{1}{2}(d-3)(d-4)+1-g$, that $C$ is contained in a unique quadric 
surface and that $\rho_{C} (n) = \rho^{S}_{d,g} (n)$ for $n \geq 
1$. Since the Rao function of a curve of degree $d$ contained in a 
quadric surface is symmetric around $\frac{d-2}{2}$, $C$ is subextremal. 
\end{proof} 
 
As a corollary, we now identify extremal and subextremal curves in 
a smooth quadric surface and in a double plane, leaving the case of 
the quadric cone and of the union of two distinct planes to the reader. 
 
\bc 
Let $C$ be an effective divisor of type $(a,b)$ on the smooth quadric surface 
$Q$, with $a \leq b$. Then  
\begin{enumerate}  
\item 
$C$ is planar if and only if  $(a,b) \in \{(0,1), (1,1) \}$;  
\item 
$C$ is extremal if and only if $(a,b) \in \{(0,2), (1,2), (2,2) \}$;  
\item 
$C$ is subextremal if and only if $(a,b) \in \{ (1,3), (2,3), (3,3) \}$. 
\end{enumerate} 
\ec

\bc 
Let $C$ be a curve in the double plane $2H$ with associated triple  
$\{Z,Y,P\}$. Then  
\begin{enumerate} 
\item 
$C$ is planar if and only if  either $y=0$, or $p=y=1$ and $z=0$; 
\item  
$C$ is extremal if and only if either $y=1$ and $p \geq 2$, or $y=p=1$ and 
$z \geq 1$, or $y=p=2$ and $z=0$;  
\item 
$C$ is  subextremal if and only if $Z$ is contained in a line and   
either $y=2$, $p \geq 3$, or $y=p=2$ and $z \geq 1$, or $y=p=3$ and $z=0$. 
\end{enumerate} 
\ec 
 
The following proposition has been proven independently by Nollet 
\cite{nolletp}; see also \cite{samir}, especially Proposition~4.15.  
 
\bp 
Let $H$ be a plane in $\pso$. 
For every $d \geq 5$ and $g \leq \frac{1}{2} (d-3)(d-4)+1$ 
(resp. $d=4$ and $g \leq 0$) the closure of the family of subextremal curves 
in $H_{d,g} (2H)$ contains an extremal curve. In particular, for $d 
\geq 4$ and $g \geq \frac{1}{2} (d-3)(d-4)+1$ the Hilbert scheme 
$H_{d,g} (\pso)$ is connected. 
\ep 
\begin{proof} 
The first statement follows from proposition~\ref{pconnect}. For $d 
\geq 4$, $g > \frac{1}{2} (d-3)(d-4)+1$ and for $(d,g)=(4,1)$ the  
Hilbert scheme $H_{d,g} (\pso)$ is irreducible \cite{extremal}, 
section 5.  
For $d \geq 5$ and $g = \frac{1}{2} (d-3)(d-4)+1$, $H_{d,g} (\pso)$    
has two irreducible components, namely the closure of the family of 
subextremal curves and the family of extremal curves \cite{samir}.  
Hence it is connected by the first statement.  
\end{proof}


\end{document}